 \documentclass[12pt,leqno]{amsart}
\newtheorem{dummy}{}[section]

\newtheorem{theorem}[dummy]{Theorem}
\newtheorem{conjecture}[dummy]{Conjecture}

\begin{document}
\bibliographystyle{plain}
\title{Derived Langlands VI: Monomial resolutions and $2$-variable L-functions }
\author{Victor P. Snaith}
\date{24 November 2020}
\maketitle
 \tableofcontents  
 
 \section{$ L(s, {\rm Ind}_{H}^{\hat{G}_{F} }(\phi) , \pi) $ for local fields}

Due to time pressure (see footnote $1$) this essay, which is part of a series consisting of a book \cite{Sn18} and 
(\cite{Sn20}, \cite{Sn20b}, \cite{newindnotes}, \cite{Sn20c}),
has a number of speculative constructions, which I do not have time to provide the rigourous details. The main constructions appear in \S1 and \S2, there is a more extensive version \cite{DLSixbigdoc} containing several other sections, being largely for reference.

Let $G$ be a (usually connected) reductive algebraic group defined over a global field $F$. Therefore 
$F$ is an algebraic number field or a function field in one variable over a finite field.

Often I shall be concerned with the points of $G$ over some local field given the completion of $F$ at a non-Archimedean prime of $F$. In this case, now writing $F$ for its local completion, suppose that $K$ is a 
(usually finite) Galois extension of $F$ and that $G$ is a quasi-split group over $F$ which splits over $K$.

I am going to use (doing my best to given page and line references for the terminology) the notation and
conventions of \cite{LLProb1970}. That is a rather tall order for the reader, should one exist, but the notation in \cite{LLProb1970} is very technical and elaborate and is explained at more length than I can manage in the time available\footnote{A deteriorating health problem over which I am not in control \cite{JMOS18}.}

Let $K/F$ be a Galois extension of non-Archimedean local fields with Galois group ${\rm Gal}(K/F)$. 

The main construction of this essay will be in the context of the following result.
\newpage
\begin{theorem}{(\cite{LLProb1970} }
\label{1.1}
\begin{em}

There is a Chevalley lattice in the Lie algebra of $G$ whose stabiliser $U_{K}$ is invariant under ${\rm Gal}(K/F)$. $U_{K}$ is self-normalising. Moreover, $G_{K} = B_{K} U_{K} $, 
$H^{1}( {\rm Gal}(K/F) ; U_{K}) = \{1 \}$ and $H^{1}( {\rm Gal}(K/F) ; B_{K} \bigcap U_{K}) = \{1 \}$. If we choose two such Chevalley lattices with stabilisers $U_{K}$ and $U'_{K}$ then $U_{K}$ and $U'_{K}$ are conjugate in $G_{K}$.
\end{em}
\end{theorem} 

The notation for Theorem \ref{1.1} is given on (\cite{LLProb1970} pp.29 final paragraph). Suppose that $K/F$ is a (not necessarily unramified for this essay) extension of local fields and $G$ is a quasi-split group over $F$ which splits over $K$. Let $B$ be a Borel subgroup of $G$ and $T$ a Cartan subgroup of $B$ both of which are defined over $F$. Let $v$ be the valuation on $K$. It is a homomorphism from $K^{*}$ whose kernel is the group of units. ${\mathcal O}_{K}^{*}$. If $t \in T_{F}$ let $v(t) \in \hat{L}$ be defined by $< \lambda , v(t) > = v( \lambda(t)) $ for all $\lambda \in L$, the group of rational characters of $T$ (\cite{LLProb1970} p,22, line -6).

The dual group $\hat{G}$ of $G$ and the Galois action on ${\rm Gal}(K/F)$ on it are defined in 
 (\cite{LLProb1970} p.22 \S2 to p.26 line 10) and, as hinted at above, the notation is quite involved but the constructions are straightforward enough. This material enables one to define 
  (\cite{LLProb1970} p.26 line 10) $\hat{G}_{F}$, which is the semi-direct product of the Galois group with $\hat{G}$.

Suppose that $\rho$ is a complex analytic representation of the semi-direct product $\hat{G}_{F}$
 (\cite{LLProb1970} p.34 line 1)\footnote{In \cite{LLProb1970} $K/F$ is taken to be an unramified extension of local fields, which is sufficient but immaterial.} and that $\pi$ is an irreducible unitary representation of $G_{F}$ on ${\mathcal H}$ whose restriction to $U_{F}$ contains the trivial representation (i.e. ${\mathcal H}^{U_{F}}$ is one-dimensional. 
 
 If $C_{c}(G_{F}, U_{F})$ is the Hecke algebra (\cite{LLProb1970} p. 30 line -7) of compactly supported functioons $f$ such that $f(gu) = f(g) = f(ug)$ for all $u \in U_{F}, g \in G_{F}$. There is a representation of $C_{c}(G_{F}, U_{F})$ on the subspace  ${\mathcal H}^{U_{F}}$ qhich gives a homomorphism $\chi$ from 
 $C_{c}(G_{F}, U_{F})$ into the ring of complex numbers  (\cite{LLProb1970} p.33 line 6). The Galois properties of this $\chi$ determine a well-defined conjugacy class in the semi-direct product $\hat{G}_{F}$ - denoted by $t \sigma_{F}$ in (\cite{LLProb1970} p.34).
  
Therefore $\rho(t \sigma_{F})$ makes sense, up to conjugate, and Langlands defines the $2$-variable $L$-function in tis case by the formula ( \cite{LLProb1970} p.34 line 4) 
\[   L(s, \rho, \pi) = \frac{1}{ {\rm det}(1 -  \rho(t \sigma_{F}) | \pi_{F}|^{s}    )}  \]
where $\pi_{F}$ generates the maximal ideal of ${\mathcal O}_{F}$. 

 Note that the "{\rm det}" in Langlands definition above is legitimatised in the sense of 
 (\cite{BAHduTM1996}  \S Theorem 3.1). A semi-simple element $a$ in the socle of a Banach algebra
 has an associated  determinant, ${\rm det}(1 - a)$ given by the formula of  (\cite{BAHduTM1996}  \S Theorem 3.1). The Banach algebra involved can be taken to be any completion containing $ \rho(t \sigma_{F}) $.

The same legitimisation, this time taking place in the $n \times n$ matrix ring with entries of the type $ \rho(t_{i,j} \sigma_{F}) $ (see next paragraph) will be required in the explanation of the  definition of 
\[ L(s, {\rm Ind}_{H}^{\hat{G}_{F} }(\phi) , \pi) \]
where $H$ is a subgroup of the semi-direct product $\hat{G}_{F}$ containing the semi-direct product of
$Z(\hat{G}_{K})$ with ${\rm Gal}(K/F)$, modulo which it is compact open, and $\phi $ is a continuous complex-values character on $H$.

Firstly I believe that the example of (\cite{DLSixbigdoc} \S4) is typical and that ${\rm Ind}_{H}^{\hat{G}_{F} }(\phi)^{U}$ is finite-dimensional. Choosing a basis $v_{1}, \ldots , v_{n}$ gives, by Langlands construction an $n \times n$  ``matrix'' $M$ of examples $t_{i,j} \sigma_{F}$ in the semi-direct product $\hat{G}_{F}$ on which the effect of making different choices is to conjugate the matrix - elementwise - by an element of  $\hat{G}_{F}$.
Let $\rho^{I_{F}}$ denote the representation of $\hat{G}_{F}$ given by the subspace of $\rho$ fixed by the semi-direct product of the Galois group with the decomposition group of $i_{F}$ of $\pi_{F}$.

Modulo the legitimisation of "{\rm det}"  my definition if the $L$-function is given by
\[   L(s,  {\rm Ind}_{H}^{\hat{G}_{F} }(\phi), \pi) = \frac{1}{ {\rm det}(1 -  \rho^{I_{F}}(M) | \pi_{F}|^{s}    )}.  \]

I expect this definition to be well-defined, to be bi-multiplicative in the variables $\rho$ and $\pi$ and whose adelic Euler product should enjoy the sort of inductivity properties which are satisfied by the Artin L-function which are recapitulated in \S\S6-10 of a more extensive version of this essay \cite{DLSixbigdoc}.

\section{Expectations of $L(s, \rho, \pi)$ via monomial resolutions}

For $\rho$ and $\pi$ admissible representations in the local field situation of \S1 suppose that
\[  \ldots  \longrightarrow  M_{i}  \longrightarrow  M_{i-1} \ldots  \longrightarrow  M_{0}   \longrightarrow  V  \longrightarrow  0 \]
is a monomial resolution (\cite{Sn18}, \cite{Sn20}  \S\S 8-10) then we are entitled, from \S1, to a $2$-variable $L$-function  $L(s,  M_{i}, \pi) $ defined as the product of the $2$-variable $L$-functions 
$ L(s,  {\rm Ind}_{H_{i,j} }^{\hat{G}_{F} }(\phi), \pi)$ such that $M_{i} = \oplus_{j} \ \underline{{\rm Ind}}_{H_{i,j} }^{\hat{G}_{F} }(\phi) $ in the monomial category.

When the representations are complex the monomial resolution is of ``finite type'' - a consequence of monomial resolutions of finite dimmensional representations of finite groups being actually finite \cite{Sn18}.
As a consequence I expect this definition of $ L(s,  \rho , \pi) $ to coincide with the definition of 
(\cite{LLProb1970} pp.29-34) and to enjoy all the analytic properties of that example.

In ( \cite{Sn20} \S4) the notion of ${\mathcal M}_{cmc, \underline{\phi}}(G)$-admissibility is introduced. It is particularly interesting in the di-p-adic situation (where the local field is $p$-adic and the representations are defined over the algebraic closure of ${\mathbb Q}_{p}$. 

In this context the monomial resolution is not necessarily of finite type, as far as I know at the moment, nevertheless I am still optimistic about the following conjecture.
\begin{conjecture}{$_{}$}
\label{2.1d}
\begin{em}

The multiplicative Euler characteristic of $2$-variable $L$-functions $L(s,  M_{i}, \pi)$ defines a well-defined and analytically well-behaved $L$-function $ L(s,  \rho , \pi) $.
\end{em}
\end{conjecture} 

When our Galois extension $K/F$ is an extension of global fields, it is explained in \cite{Sn18} how to use the 
adelic Tensor Product Theorem (\cite{Sn18} Theorem 1.21) to define adelic monomial resolutions of automorphic representations. In this situation there is a simple reduction to the case in which the finite Galois group is soluble, which is based on the fact that the subgroups involved in the monomial resolution of any finite-dimensional of a finite group over any algebraically closed field of characteristic zero may be taken to be $M$-groups (\cite{Sn94} Proposition 2.1.17 p.30). 

Recall that a finite group $G$ is nilpotent if and only if it has a lower central series
\[ \{ 1 \} = Z_{0} \lhd  Z_{1} = Z(G) \lhd Z_{2} \lhd \ldots \lhd Z_{n} =G \]
exists such that $Z_{i+1}/Z_{i} = Z(G/Z_{i})  $ for all $i$. In particular nilpotent groups are $M$-groups - each irreducible representation is induced from a $1$-dimensional character of a subgroup. Since nilpotent groups are the product of their Sylow $p$-subgroups a subgroup of a nilpotent group is again nilpotent - but not so for M-groups.

An M-group is soluble \cite{MI1994}. The derived series of $G$ is 
\[ G^{0} = G, G^{1} = [G^{0}, G^{0}] = [G,G], {\rm the \ commutator \ subgroup \ of \ G}  \]
and $ G^{n} = [G^{n-1}, G^{n-1 }] $ for $n= 1,2,3,\ldots $. A soluble group $G$ is one for which $G^{n} = \{ 1 \}$ for some $n$. Subgroups of soluble groups are soluble since $J < G$ implies $[J,J] < [G,G]$.
We begin by writing the trivial one-dimensional representation of the Galois group as a sum of iduced monomial representations ${\rm Ind}_{H}^{G}(\phi)$ with $H$ an $M$-subgroup of $G$. We can inflate this relation to the semi-direct prodct and multiply the monomial resolution, as a complex of representation, with this formula for $1$.

Recall also that the product in $R_{+}(G)$ is given by a double coset formula (\cite{Sn94} p.68, Exercise 2.5.7) 
\[ \begin{array}{l}
(K, \phi)^{G} \cdot (H, \psi)^{G} \\
\\
=  \sum_{w \in K  \backslash G /H} \ ( w^{-1} Kw \bigcap  H, w^{*}(\phi) \psi )^{G} .
\end{array} \]
Suppose, abbreviating ${\rm Gal}(K/F)$ to ${\rm Gal}$, we have a Galois semi-direct product ${\rm Gal} \propto G$ with $(K, \phi)^{{\rm Gal} \propto G} \in R_{+}({\rm Gal} \propto G) $ and if $(H, \psi)^{{\rm Gal}}$ and $\lambda : {\rm Gal} \propto G \longrightarrow   {\rm Gal}$ is the projection then $( \lambda^{-1}(H) , \psi \lambda)^{{\rm Gal}  \propto G}$. Therefore 
\[ \begin{array}{l}
(K, \phi)^{{\rm Gal} \propto G} \cdot ( \lambda^{-1}(H) , \psi \lambda)^{{\rm Gal}  \propto G} \\
\\
=  \sum_{w \in K  \backslash {\rm Gal} \propto G /H \propto G} \ ( w^{-1} Kw \bigcap  H, w^{*}(\phi) \psi )^{G} .
\end{array} \]
If $w = (z,g) \in  {\rm Gal} \propto G$ and $(y, g_{1} ) \in H \propto G$ then $ (z,g)(y, g_{1} ) = ( zy, g z(g_{1}))$so we may take $w = (z,1) \in  K \bigcap {\rm Gal} \backslash {\rm Gal} / H$ and $w^{-1} K w \bigcap
 H \propto G \subset H \propto G$. 
 
This procedure reduced the Galois semi-direct products to ones involving solution Galois groups.


\begin{thebibliography}{}

  \bibitem{BAHduTM1996}  B. Aupetit and H. du T. Mouton: Trace and determinant in Banach algebras; Studia Mathematica 121 (2) 114-136 (1996).

  \bibitem{LLProb1970}  R.P. Langlands: Problems in the theory of automorphic forms; {\em Lectures in modern analysis and its applications III}, Springer LNMath \#170 (1970) 18-61.
  
   \bibitem{MI1994} M. I. Isaacs: {\em Characters of finite groups}; Dover (1994) ISBN  978-0-486-68014-9.
   
   \bibitem{JMOS18} J.M. O'Sullivan and C.N. Harrison: Myelofibrosis: Clinicopathologic Features, Prognosis and Management; Clinical Advances in Haematology and Oncology {\bf 16} (2) February 2018.1.

  
  \bibitem{Sn94}  V.P. Snaith:  {\em Explicit Brauer Induction (with applications to algebra and number theory)}; Cambridge studies in advanced mathematics \#40,Cambridge University Press (1994).

\bibitem{Sn18} V.P. Snaith:  {\em Derived Langlands}; World Scientific (2018).

\bibitem{Sn20} V.P. Snaith:  Derived Langlands II: HyperHecke algebras, monocentric relatons and ${\mathcal M}_{cmc, \underline{\phi}}(G)$-admissibility; arXiv:3100675 [math.RT] 24 Mar 2020.

\bibitem{Sn20b} V.P. Snaith: Derived Langlands III: PSH algebras and their numerical invariants; arXiv:3223311 [math.RT] 12 June 2020.

\bibitem{newindnotes}  Derived Langlands IV: Notes on ${\mathcal M}_{c}(G)$-induced representations; arXiv:2008.06325v1 [math.RT] 14 Aug 2020.

\bibitem{Sn20c} V.P. Snaith: Derived Langlands V: The Hopflike properties of the hyperHecke algebra; preprint on University of Sheffield homepage (27 May 2020).

\bibitem{DLSixbigdoc}  V. P. Snaith: Derived Langlands VI: Monomial resolutions and $2$-variable L-functions
(extensive version);  http://victor-snaith.staff.shef.ac.uk/   (November 2020).

\end{thebibliography}
 \end{document}